\documentclass[]{theclass}
\begin{document}
\begin{frontmatter}

\titledata{An updated survey on $2$-Factors of Regular Graphs}{}           

\authordatatwo{Domenico Labbate}{domenico.labbate@unibas.it}{}{Federico Romaniello}{federico.romaniello@unibas.it}{}
{Dipartimento di Matematica, Informatica ed Economia\\ Universit\`{a} degli Studi della Basilicata, Potenza, Italy}


\begin{abstract}
A 2-factor of a graph $G$ is a 2-regular spanning subgraph of $G$. We present a survey summarising results on the structure of 2-factors in regular graphs, as achieved by various researchers in recent years.
\end{abstract}

\end{frontmatter}

\section{Introduction}

All graphs considered in this survey are finite and simple (without loops or multiple edges). We shall use the term \textit{multigraph} when multiple edges are permitted. For definitions and notations not explicitly stated the reader may refer to Bondy and Murty's book \textit{Graph Theory}\cite{BM}.

Several authors have considered the number of Hamiltonian circuits in $k$--regular graphs and there are interesting and beautiful results and conjectures in the
literature. In particular, C.A.B. Smith (1940, cf. Tutte \cite{T})
proved that each edge of a $3$--regular multigraph lies in an
even number of
Hamiltonian circuits. This result was extended to multigraphs in which
each vertex has an odd
degree by Thomason \cite{Tho}.

A multigraph with
exactly one Hamiltonian circuit is said to be {\em Uniquely Hamiltonian}.
Thomason's result implies that there are no regular
Uniquely Hamiltonian multigraphs of odd degree.
In 1975, Sheehan \cite{S} posed the following famous conjecture:

\begin{conj}\cite{S}\label{sheehan con}
There are no Uniquely Hamiltonian $k$--regular
graphs for all integers $k \geq 3$.
\end{conj}

It is well known that it is enough to prove it for $k=4$. This conjecture
has been verified by Thomassen for bipartite graphs, \cite{Th1}
(under the weaker hypothesis that $G$ has minimum degree $3$) and for $k$--regular graphs when $k \geq 300$, \cite{Th3}.
Ghandehari and Hatami have improved this value for  $k \geq 48$ \cite{Ghand} and, recently, by Haxell, Seamone, and Verstraete \cite{HaxSea} for  $k > 22$.

In this context, several recent papers addressed the problem of characterising
families of graphs (particularly regular) that have certain conditions imposed on their 2--factors.
This survey presents the main results obtained in the recent years.
We will also discuss the connections of these problems with the particular class of {\em odd 2--factored snarks}.

\section{Preliminaries}

An ${\it r}${\em --factor} of a graph $G$ is an $r$--regular
spanning subgraph of $G$. Thus, a $2${\em --factor} of a graph $G$ is a
$2$--regular spanning subgraph of $G$, while a 1-factor is also called \emph{perfect matching} since it is a matching that covers all the vertices. 
A {\it 1--factorization} of $G$ is a
partition of the edge set of $G$ into edge-disjoint $1$--factors.

Let $G$ be a bipartite graph with bipartition
$(X,Y)$ such that $|X|=|Y|$, and $A$ be its adjacency
matrix. In general $0 \leq |det(A)| \leq per(A)$.
We say that $G$ is {\em  det--extremal} if $|det(A)|=per(A)$.
Let $X=\{x_1,x_2,\ldots,x_n\}$ and $Y=\{y_1,y_2,\ldots,y_n\}$
be the bipartition
of $G$. For $F$ a $1$--factor of $G$ we define the {\em sign}
of $F$, $sgn(F)$, to be the sign of the permutation
of $\{1,2,\ldots,n\}$ corresponding to $F$. (Thus $G$ is det--extremal
if and only if all $1$--factors of $G$ have the same sign.)
The following elementary result is a special case of
\cite[Lemma 8.3.1]{LP}.
\begin{lemma}\label{sgn}
Let $F_1,F_2$ be 1-factors in a bipartite graph $G$
and $t$ be the number of circuits in
$F_1\cup F_2$ of length congruent to zero modulo four.
Then
$sgn(F_1)sgn(F_2)=(-1)^t$.
\end{lemma}

\begin{figure}[h]
\begin{center}
\includegraphics[scale=0.18]{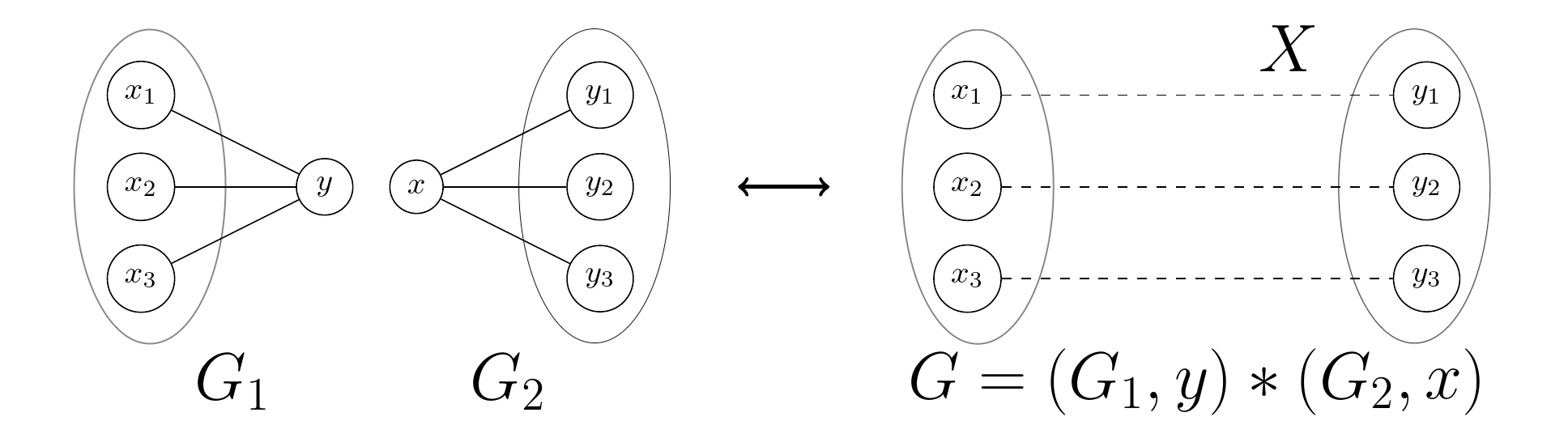}
\caption{From left to right: The Star Product between the two graphs $G_1$ and $G_2$. From right to left: The 3-edge-reduction on $X$ in the graph $G$.}
\label{starp}
\end{center}
\end{figure}

Before proceeding, we recall a standard operation on graphs that will be recurrent in this survey. Let $G_1$ and $G_2$ be two graphs each containing a vertex of degree $3$, say $y\in V(G_1)$ and $x\in V(G_2)$.
Let $x_1,x_2,x_3$ be the neighbours of $y$ in $G_1$ and
$y_1,y_2,y_3$ be the neighbours of $x$ in $G_2$.
We say that the graph $G=(G_1-y)\cup (G_2-x)\cup\{x_1y_1,x_2y_2, x_3y_3\}$ is a {\em star product} of $G_1$ and $G_2$ and
write $G=(G_1,y)*(G_2,x)$. We remark that if $G_1$ and $G_2$ are bridgeless and cubic, then the graph obtained is also bridgeless and cubic. For a bridgeless cubic graph $G$ having a 3-edge-cut $X$, in the opposite direction, it is possible to define a \emph{3-edge-reduction} on $X$ as the graph operation on $G$ which creates two new bridgeless cubic graphs by adding a new vertex to each of the components of $G - X$ and joining it to the degree two vertices in the respective component. Both operations are shown in Figure \ref{starp}.


The {\em Heawood graph} $H_0$ is the bipartite graph associated with the
point/line incidence matrix of the Fano plane $PG(2,2)$.
Let ${\cal{SP}}(H_0)$ be the class of graphs obtained from the Heawood graph by
repeated star products.

These graphs were used by McCuaig in \cite{MCC2} to characterise
the $3$--connected cubic det-extremal bipartite graphs:

        \begin{theorem}\label{3-cont}\cite{MCC2}
A 3-connected cubic bipartite graph
is det-extremal if and only if it belongs to ${\cal{SP}}(H_0)$.
        \end{theorem}

\vspace{0.125in}
\noindent
{\bf Note} (i) Theorem \ref{3-cont} has been improved for connectivity $2$ graphs by Funk, Jackson, Labbate and Sheehan in \cite{FJLS1};

(ii) Bipartite graphs $G$ with the more general property that
some of the entries in the adjacency matrix $A$ of $G$ can be changed from
$1$ to $-1$ in such a way that the resulting matrix $A^{\ast}$ satisfies
$per(A)=det(A^{\ast})$ have been characterised in \cite{Li,MCC,RST}.

\section{2--factor Hamiltonian graphs}\label{2factham}

A graph admitting with a $2$--factor is said to be {\em $2$--factor
Hamiltonian} if all its $2$--factors are Hamiltonian cycles.
Examples of such graphs are $K_4$, $K_5$, $K_{3,3}$, the Heawood graph $H_0$, and the cubic graph of girth five
obtained from a 9-circuit by adding three vertices, each joined to three
vertices of the 9-circuit, also known as the `Triplex graph' of Robertson,
Seymour and Thomas.

The following property was stated without proof in \cite{FJLS2}, and it is important for approaching a characterization of this family of graphs.

        \begin{prop}\label{star}
If $G$ is a bipartite graph represented
as a star product $G=(G_1,y)*(G_2,x)$ and every pair of edges in the 3-edge cut of $G$ belong to a 2--factor. Then $G$ is 2--factor Hamiltonian if
and only if $G_1$ and $G_2$ are 2--factor Hamiltonian.
        \end{prop}

We have stated here a slightly different version of the above Proposition \ref{star}, In fact, we have added here the weaker hypothesis that every pair of edges in the 3-edge cut of $G$ belong to a 2-factor, which was tacitly assumed in \cite{FJLS2}. 

The bipartite hypothesis in Proposition \ref{star} is needed, as the star products of
non--bipartite $2$--factor Hamiltonian graphs is not necessarily $2$--factor Hamiltonian. It happens, for example, when considering $K_4 \ast K_4$.

In addition to this, it was pointed out by M. Gorsky and T. Johanni in a private communication \cite{gtpc} that without the hypothesis on the 3-edge cut of $G$ it is possible to obtain 2-factor Hamiltonian graphs as a star product of two graphs that are not 2-factor Hamiltonian, as shown in Figure \ref{fig:count2fh}.\\
For these reasons, a proof of the above Proposition \ref{star} is given below.
\begin{proof}
We start by noticing that since we are assuming that $G$ is bipartite, then so are $G_1$ and $G_2$, by the properties of the star product. Let $X$ be the 3-edge cut of the bipartite graph $G$.\\
($\Rightarrow$) It follows immediately that the 2--factors in $G_1$ and $G_2$ arises from 2--factors of $G$ by contracting the edges of the 2--factor in $X$. Since $G$ is 2--factor Hamiltonian, then $G_1$ and $G_2$ must be at least Hamiltonian graphs. Suppose now, by contradiction that $G_1$ or $G_2$ are not 2--factor Hamiltonian, say $G_1$. Hence, there exists a 2--factor $F_1$ of $G_1$ made of at least two disjoint cycles. Let $C$ be the cycle in $F_1$ containing the vertex $y$ on which the star product is operating. It follows that the two edges of $C$ incident to $y$ correspond to two edges of $X$. Based on our assumption, these two edges are part of a 2--factor of $G$, enabling us to combine $C$ with the Hamiltonian cycle of $G_2$ to form a 2-factor of $G$ with at least two components. A contradiction.\\
($\Leftarrow$) Suppose now that $G_1$ and $G_2$ are 2--factor Hamiltonian and, by contradiction suppose that $G$ is not 2--factor Hamiltonian. Note that by hypothesis $G$ admits 2--factors containing edges of $X$. Let $F$ be a 2--factor of $G$ which is not a Hamiltonian cycle. It follows that $|F \cap X|$ must be equal to two or to zero. In the former case, it follows immediately that $G_1$ or $G_2$ are not 2--factor Hamiltonian, which is a contradiction. Suppose now that $|F \cap X| = 0$. Since $G$ is a bipartite graph containing a 2--factor it has an even number of vertices, and $G - X$ has exactly two components. It follows that both $G_1$ and $G_2$ are of odd order, but since they are Hamiltonian they contain an odd cycle, which is a contradiction.
\end{proof}

\begin{figure}[h!]
\centering
\includegraphics[scale=0.15]{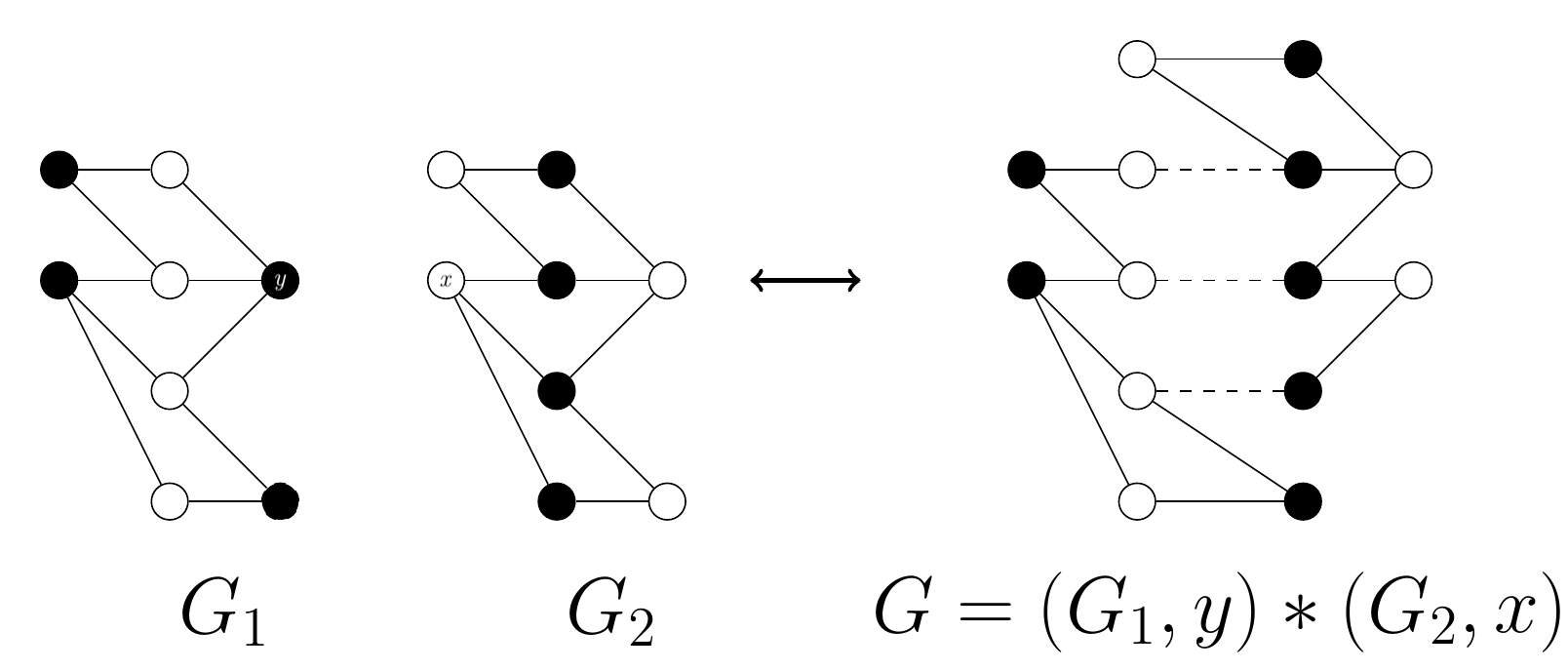}
\caption{The star product $(G_1,y) \ast (G_2,x)$ between the two non-2-factor Hamiltonian graphs $G_1$ and $G_2$ is 2-factor Hamiltonian. However, not every pair of the dashed edges of $X$ belong to a 2--factor, as shown in \cite{gtpc}.}
\label{fig:count2fh}
\end{figure}

\

Indeed, it is worth remarking the additional hypothesis on the 3-edge cut of $G$ is not needed when dealing with bipartite cubic graphs as it is always satisfied, and by using Proposition \ref{star}, Funk, Jackson, Labbate, and Sheehan constructed an infinite family of 2--factor Hamiltonian
cubic bipartite graphs by taking iterated  star products of
$K_{3,3}$ and $H_0$ \cite{FJLS2}.
They conjecture that these are the only non-trivial
2--factor Hamiltonian regular bipartite graphs.

        \begin{conj}\label{conj1}\cite{FJLS2}
Let $G$ be a $2$--factor Hamiltonian
$k$-regular bipartite graph. Then either $k=2$ and $G$ is a circuit or
$k=3$ and $G$ can be obtained from $K_{3,3}$ and $H_0$ by repeated star
products.
        \end{conj}

We remark here that Conjecture \ref{conj1} is still open, and a positive answer to it will allow to completely characterise the family of 2--factor Hamiltonian regular bipartite graphs.
In the 1980s Sheehan posed the following Conjecture:

\begin{conj}\label{conj2}
There are no 2--factor Hamiltonian k--regular bipartite graphs for all integers $k\geq 4$.
\end{conj}

The following properties have been proved by Labbate \cite{L1,L2} for an equivalent family of cubic graphs (cf. subsection \ref{m1f-sect}), and then by Funk, Jackson, Labbate and Sheehan \cite{FJLS2} for 2--factor Hamiltonian graphs:

\begin{lemma}\label{m1f}\cite{L2,L1,FJLS2}
Let $G$ be a $2$--factor Hamiltonian cubic bipartite graph.
Then $G$ is $3$--connected and $|V(G)| \equiv 2 \pmod{4}$.
\end{lemma}

A graph $H$ is {\em `maximally'} 2--factor Hamiltonian if the multigraph $G$ obtained by
adding an edge $e$ with endvertices $u,v$ to $H$ has a disconnected $2$--factor containing $e$.

\begin{lemma}\label{maximal}\cite[Lemma 3.4 (a)(i)]{FJLS2}
Graphs obtained by taking star products of $H_0$ are maximally 2--factor Hamiltonian.
\end{lemma}

Funk, Jackson, Labbate and Sheehan in \cite{FJLS2} proved Conjecture \ref{conj2} applying Lemmas \ref{sgn}, \ref{m1f}, \ref{maximal} and Theorem \ref{3-cont}:

        \begin{theorem}\label{existence}\cite{FJLS2}
Let $G$ be a 2--factor Hamiltonian $k$--regular bipartite graph.
Then $k\leq 3$.
        \end{theorem}

Theorem \ref{existence} has inspired further results by Faudree, Gould, and Jacobsen \cite{FGJ} that determined the maximum number of edges in both
$2$--factor Hamiltonian graphs and $2$--factor Hamiltonian
bipartite graphs. In particular, they proved the following theorems:

    \begin{theorem}\label{FGJbip}\cite{FGJ}
If $G$ is a bipartite $2$--factor Hamiltonian graph of order $n$ then
$$
|E(G)| \, \leq \, \left\{
\begin{array}{lcc}
n^2/8+n/2     & if & n \equiv 0 \mod 4, \\
n^2/8+n/2+1/2 & if & n \equiv 2 \mod 4,
\end{array}
\right.
$$
and the bound is sharp.
    \end{theorem}

    \begin{theorem}\label{FGJnonbip}\cite{FGJ}
If $G$ is a $2$--factor Hamiltonian graph of order $n$ then
$$
|E(G)| \, \leq \, \lceil n^2/4+n/4  \rceil
$$
and the bound is sharp for all $n\geq 6$.
    \end{theorem}

In addition, Diwan \cite{Di} has shown that

    \begin{theorem}\label{diwanplanar}
$K_4$  is the only $3$--regular $2$--factor Hamiltonian planar
graph.
    \end{theorem}

Conjecture \ref{conj1} has been partially solved in terms of minimally 1--factorable cubic bipartite graphs as we explain in the following subsection.

\subsection{Minimally 1-factorable graphs}\label{m1f-sect}

Let $G$ be a $k$--regular bipartite graph.
We say that $G$ is {\em minimally 1--factorable}
if every 1--factor of $G$ is contained in a unique
1--factorization of $G$.

The results cited above by Funk, Jackson, Labbate, and Sheehan were inspired by results on minimally 1-factorable graphs
obtained in \cite{FL,L1,L2,L3}. It can be seen that:

\begin{prop}\label{equiv}\cite{FJLS2}
If $G$ is minimally $1$--factorable then $G$ is 2--factor Hamiltonian.
If $k=2,3$, then $G$ is minimally $1$--factorable if and only if
$G$ is 2--factor Hamiltonian.
\end{prop}

Theorem \ref{existence} extends the result of \cite{FL} that
minimally $1$--factorable $k$--regular bipartite graphs
exist only when $k\leq 3$.

Furthermore, Labbate in \cite{L3} proved the following characterization:

\begin{theorem}\label{char}\cite{L3}
Let $G$ be a minimally 1--factorable $k$--regular bipartite graph of girth $4$. Then either $k=2$ and $G$ is a circuit or
$k=3$ and $G$ can be obtained from $K_{3,3}$ by repeated star
products.
\end{theorem}

Hence, it follows from results in \cite{L2} that a smallest counterexample to Conjecture \ref{conj1}
is cubic and cyclically 4-edge connected, and from Theorem \ref{char}
that it has girth at least six. Thus, to prove
the conjecture, it would suffice to show that the Heawood graph is the only
2-factor Hamiltonian cyclically 4--edge connected cubic bipartite
graph of girth at least six.

This seems a challenging task to achieve at least with the techniques used so far.
In \cite{ALS}, partial results were obtained by using {\em irreducible Levi graphs} (cf. Section \ref{irredlevi} and Theorem \ref{MainThm}).

\subsection{Perfect Matching Hamiltonian Graphs}
A graph $G$ admitting a $1$--factor is said to have the \textit{Perfect-Matching-Hamiltonian property} (for short, the PMH-property) if every $1$--factor $M$ of $G$ can be extended to a Hamiltonian cycle of $G$, that is, there exists a $1$--factor $N$ of $G$ such that $M\cup N$ induces a Hamiltonian cycle of $G$. This problem was first introduced by Las Vergnas \cite{LV} and H\"{a}ggkvist \cite{Haggkvist} in the 1970s, and in recent years this concept was studied with particular focus on cubic graphs, the reader may find more details in \cite{papillon} and \cite{AAAHST} .  For simplicity, a graph admitting the PMH-property is said to be PMH or a PMH-graph. The class of $2$-factor Hamiltonian graphs are of course PMH-graphs. Recently, in \cite{btxbtw}, in an attempt to look at Conjecture \ref{conj1} from a different point of view, the following statement was proved:
\begin{prop}\label{betw}\cite{btxbtw}
Let $G$ be a cubic PMH-graph (not necessarily bipartite). The graph $G$ is not 2-factor Hamiltonian if and only if it admits a $1$--factor which can be extended to a Hamiltonian cycle in exactly one way.   
\end{prop}
Proposition \ref{betw} suggests another way how one can look at Conjecture \ref{conj1}. Indeed, a smallest counterexample to this  conjecture can  be searched  for  in  the  class  of  bipartite  cubic PMH-graphs, and hence the Conjecture \ref{conj1} can be equivalently restated  in  terms  of  a  strictly  weaker property than 2-factor Hamiltonicity:  the PMH-property.
\begin{conj}\cite{btxbtw}
Every bipartite cyclically 4-edge-connected cubic PMH-graph with girth at least 6, except the Heawood graph, admits a $1$--factor which can be extended to a Hamiltonian cycle in exactly one way.   
\end{conj}

\section{2--factor isomorphic graphs}\label{2factiso}

The family of 2--factor Hamiltonian k--regular graphs can be extended to the family of
connected $k$--regular graphs with the more general property that all their
$2$--factors are isomorphic, i.e. the family of {\em 2--factor isomorphic} $k$--regular bipartite graph.

Examples of such graphs are given by all the 2--factor Hamiltonian and the Petersen graph (which is 2--factor isomorphic since it has all its 2--factors of length $(5,5)$ but it is not 2-factor Hamiltonian). Note that star product preserves also 2--factor isomorphic regular graphs.

In \cite{AFJLS} Aldred, Funk, Jackson, Labbate and Sheehan proved the following existence theorem

    \begin{theorem}\label{existencebuk}\cite{AFJLS}
Let $G$ be a 2--factor isomorphic $k$--regular bipartite graph.
Then $k\leq 3$.
    \end{theorem}

They also conjecture that the family of 2--factor isomorphic and the one of 2--factor Hamiltonian
$k$--regular bipartite graphs are, in fact, the same.

\begin{conj}\label{conjequiv}\cite{AFJLS} Let $G$ be a connected $k$-regular bipartite graph.
Then $G$ is $2$--factor isomorphic if and only if $G$ is 2--factor Hamiltonian.
\end{conj}

Abreu, Diwan, Jackson, Labbate, and Sheehan proved in \cite{ADJLS} that Conjecture \ref{conjequiv} is false by applying the following construction:

\begin{prop}\label{counterexample}\cite{ADJLS}
Let $G_i$ be a 2-factor Hamiltonian cubic bipartite graph
with $k$ vertices and
$e_i=u_iv_i\in E(G_i)$ for $i=1,2,3$.
Let $G$
be the graph obtained from the disjoint union of
the graphs $G_i-e_i$ by adding two
new vertices $w$ and $z$ and new edges
 $wu_i$ and $zv_i$ for $i=1,2,3$.
Then $G$ is a
non-Hamiltonian connected 2-factor isomorphic cubic bipartite
graph of edge-connectivity two.
\end{prop}

Given a set $\{G_1,G_2,\ldots,G_k\}$ of 3--edge--connected cubic bipartite graphs, let denote by ${\cal{SP}}(G_1,G_2,\ldots,G_k)$ the set of
cubic bipartite graphs which can be obtained from
$G_1,G_2,\ldots,G_k$ by repeated star products.
In Section \ref{2factham} we have seen that
it was shown in \cite{FJLS2} that all graphs in ${\cal
SP}(K_{3,3},H_0)$ are 2-factor Hamiltonian. Thus we may apply
Proposition \ref{counterexample} by taking $G_1=G_2=G_3$ to be any
graph in ${\cal SP}(K_{3,3},H_0)$ to obtain an infinite family of
2--edge--connected non--Hamiltonian 2--factor isomorphic cubic
bipartite graphs. This family gives counterexamples to the
Conjecture \ref{conjequiv}. Note, however, that Conjecture \ref{conjequiv} can be modified as follows:

    \begin{conj}\label{conjequivmodif}\cite{ADJLS}
Let $G$ be a 3--edge--connected 2--factor isomorphic cubic bipartite graph. Then $G$ is a
2--factor Hamiltonian cubic bipartite graph.
    \end{conj}

Recall that a \textit{digraph} is a graph in which the edges have a direction (and they are now ordered pairs of vertices). In \cite{AAFJLS, AAFJLS2} Abreu, Aldred, Funk, Jackson, and Sheehan also proved existence theorems for the digraphs and non--bipartite graphs case, as shown below.

For $v$ a vertex of a digraph $D$, let $d^+(v)$ and $d^-(v)$
denote the out--degree and in--degree of $v$.
We say that D is $k$--diregular if for all vertices $v$ of $G$, we have $d^+(v)=d^-(v)=k$.

\begin{theorem}\label{theorem3.1}\cite{AAFJLS, AAFJLS2}
Let $D$ be a digraph with $n$ vertices and
$X$ be a directed $2$-factor of $D$.
Suppose that either\\
(a) $d^+(v)\geq \lfloor\log_2 n\rfloor+2$ for all $v\in V(D)$, or\\
(b) $d^+(v)=d^-(v)\geq 4$ for all $v\in V(D)$. \\
Then $D$ has a
directed $2$-factor $Y$ with $Y\not\cong X$.
\end{theorem}

    \begin{cor}\label{eqdigraph}\cite{AAFJLS}
Let $G$ be a k--diregular directed graph. Then $k \leq 3$.
    \end{cor}

\begin{theorem}\label{theorem4.1}\cite{AAFJLS, AAFJLS2}
Let $G$ be a graph with $n$ vertices and
$X$ be a $2$-factor of $G$.
Suppose that either\\
(a) $d(v)\geq 2(\lfloor\log_2 n\rfloor+2)$ for all $v\in V(G)$, or\\
(b) $G$ is a $2k$-regular graph for some $k\geq 4$. \\
Then $G$ has a $2$-factor $Y$ with $Y\not\cong X$.
\end{theorem}

They have also posed the following open problems and conjecture:

\begin{quest}\cite{AAFJLS2} Do there exist $2$--factor isomorphic bipartite graphs
of arbitrarily large minimum degree?
\end{quest}

\begin{quest}\cite{AAFJLS2} Do there exist $2$--factor isomorphic regular graphs
of arbitrarily large degree?
\end{quest}

        \begin{conj}\label{hukemptyset}\cite{AAFJLS}
The graph $K_{5}$ is the only 2--factor Hamiltonian 4--regular non--bipartite graph.
        \end{conj}

\section{Pseudo 2--factor isomorphic graphs}

In \cite{ADJLS} Abreu, Diwan, Jackson, Labbate, and Sheehan  extended  the above-mentioned results on regular
$2$--factor isomorphic bipartite graphs to the more general family
of {\em pseudo $2$--factor isomorphic graphs} i.e. graphs $G$ with the
property that the parity of the number of circuits in a
$2$--factor is the same for all $2$--factors of $G$.

Examples of such graphs are given by all the 2--factor isomorphic regular graphs and the Pappus graph (i.e. the point/line incidence graph of the Pappus configuration). The family of pseudo $2$--factor isomorphic is wider than the one of $2$--factor isomorphic regular bipartite graphs:

    \begin{prop}\label{notequal}\cite{ADJLS}
The Pappus graph $P_0$ is pseudo $2$--factor isomorphic but not
$2$--factor isomorphic.
     \end{prop}

In \cite{ADJLS} Abreu, Diwan, Jackson, Labbate and Sheehan proved the following existence theorem:

    \begin{theorem}\label{Main}\cite{ADJLS}
Let $G$ be a pseudo $2$--factor
isomorphic $k$--regular bipartite graph. Then $k\in \{2,3\}.$
    \end{theorem}

They have also shown that there are
no planar pseudo 2--factor isomorphic
cubic bipartite graphs.

\begin{theorem}\label{nonplanar}\cite{ADJLS}
Let $G$ be a
pseudo 2--factor isomorphic cubic bipartite graph.
Then $G$ is non-planar.
\end{theorem}

Star products preserve also the property of being
pseudo $2$--factor isomorphic in the family of cubic bipartite
graphs.

    \begin{lemma}\label{starpseudo}\cite{ADJLS}
Let $G$ be a star product of two pseudo $2$--factor isomorphic
cubic bipartite graphs $G_1$ and $G_2$. Then $G$ is also pseudo
$2$--factor isomorphic.
\end{lemma}

Thus $K_{3,3}$, $H_0$ and $P_0$ can be used to construct
an infinite family of 3--edge--connected pseudo $2$--factor
isomorphic cubic bipartite graphs.

Lemma \ref{starpseudo} implies that all graphs in ${\cal{SP}}(K_{3,3},H_0,P_0)$ are pseudo $2$--factor isomorphic. In \cite{ADJLS} Abreu, Diwan, Jackson, Labbate and Sheehan
conjectured that these are the only 3--edge--connected pseudo
$2$--factor isomorphic cubic bipartite graphs.

    \begin{conj}\label{quest1}\cite{ADJLS}
Let $G$ be a 3--edge--connected cubic bipartite graph.
Then $G$ is  pseudo $2$--factor
isomorphic
if and only if $G$ belongs to ${\cal SP}(K_{3,3},H_0,P_0)$.
    \end{conj}

Recall that McCuaig \cite{MCC2} has shown that a
3--edge--connected cubic bipartite graph $G$ is det-extremal if and only if
$G\in {\cal SP}(H_0)$.

Let $G$ be a graph and $E_1$ be an edge-cut of $G$. We say that
$E_1$ is a {\em non-trivial edge-cut} if all components of $G-E_1$
have at least two vertices. The graph $G$ is {\em essentially
4--edge--connected} if $G$ is 3--edge--connected and has no
non-trivial 3--edge--cuts.
Let $G$ be a cubic bipartite graph with bipartition $(X,Y)$ and $K$ be a
non-trivial 3--edge--cut of $G$. Let $H_1,H_2$ be the components of
$G-K$. We have seen that $G$ can be expressed as a star product
$G=(G_1,y_K)*( G_2,x_K)$ where $G_1-y_K=H_1$ and $G_2-x_K=H_2$.
We say that $y_K$, respectively $x_K$, is the {\em marker vertex} of
$G_1$, respectively $G_2$, {\em corresponding to the cut $K$}.
Each non-trivial 3--edge--cut of $G$ distinct from $K$ is a
non-trivial 3--edge--cut of $G_1$ or $G_2$, and vice versa.
If $G_i$ is not essentially 4--edge--connected for $i=1,2$, then we may
reduce $G_i$ along another non-trivial 3-edge-cut. We can continue this
process until all the graphs we obtain are essentially 4--edge--connected.
We call these resulting graphs the {\em constituents} of $G$.
It is easy to see that the constituents of $G$ are unique i.e. they are
independent of the order we choose to reduce the non-trivial
3--edge--cuts of $G$.

It is also easy to see that Conjecture
\ref{quest1} holds if and only if Conjectures \ref{quest2} and
\ref{quest3} below are both valid.

\begin{conj}\label{quest2}\cite{ADJLS}
Let $G$ be an essentially $4$--edge--connected pseudo $2$--factor
isomorphic
cubic bipartite graph.
Then $G\in \{K_{3,3},H_0,P_0\}$.
\end{conj}

\begin{conj}\label{quest3}\cite{ADJLS}
Let $G$ be a $3$--edge--connected
pseudo $2$--factor isomorphic cubic bipartite graph and suppose that
$G=G_1*G_2$.
Then $G_1$ and $G_2$ are both
pseudo $2$--factor
isomorphic.
\end{conj}

In \cite{ADJLS} Abreu, Diwan, Jackson, Labbate, and Sheehan obtained partial results on Conjectures \ref{quest2} and \ref{quest3}
as follows:

\begin{theorem}\label{ess4ec}\cite{ADJLS}
Let $G$ be an essentially $4$--edge--connected pseudo 2--factor isomorphic
cubic bipartite graph. Suppose $G$ contains a $4$--circuit. Then $G=K_{3,3}$.
\end{theorem}

They used Theorem \ref{ess4ec} to deduce some evidence in favour of Conjecture \ref{quest1}.

\begin{theorem}\label{constit}\cite{ADJLS}
Let $G$ be a $3$--edge-connected pseudo 2--factor isomorphic bipartite graph.
Suppose $G$ contains a 4--cycle $C$. Then $C$ is contained in a constituent of $G$
which is isomorphic to $K_{3,3}$.
\end{theorem}

Note that Theorem \ref{constit} generalises Theorem \ref{char}
obtained by Labbate in \cite{L3} for minimally 1--factorable bipartite cubic graphs (or equivalently 2--factor Hamiltonian cubic bipartite graphs) to the family of pseudo 2--factor isomorphic bipartite graphs.
Furthermore, Theorem \ref{constit} leaves the characterization of pseudo 2--factor isomorphic bipartite graph open for girth $\ge 6$.

Abreu, Labbate and Sheehan \cite{ALS1} gave a partial solution to this open case in terms of irreducible configuration of Levi graphs as described in the next subsection.

\subsection{Irreducible pseudo 2--factor isomorphic cubic bipartite graphs}\label{irredlevi}

An incidence structure is {\em linear} if two different points are
incident with at most one line. A  {\em  symmetric configuration} $n_k$ (or {\em $n_k$ configuration})
is a linear incidence structure consisting of $n$ points
and $n$ lines such that each point and line is respectively
incident with $k$ lines and points.
Let $\cal C$ be a symmetric configuration $n_k$, its {\em Levi graph} $G({\cal C})$
is a $k$--regular bipartite graph whose vertex set are the points and the lines of $\cal C$ and
there is an edge between a point and a line in the graph if and only if they are incident in $\cal C$.
We will indistinctly refer to Levi graphs of configurations as their {\em incidence graphs}.

It follows from Theorem \ref{constit} that an essentially
$4$--edge--connected pseudo $2$--factor isomorphic cubic bipartite
graph of girth greater than or equal to $6$ is the Levi graph of
a symmetric configuration $n_3$.

In $1886$ V. Martinetti \cite{VM}
characterised symmetric configurations $n_3$, showing that they can be
obtained from an infinite set of so called {\em irreducible}
configurations, of which he gave a list.
Recently, Boben proved that Martinetti's list of irreducible
configurations was incomplete and completed it \cite{B}.
Boben's list of irreducible configurations was obtained by characterising
their Levi graphs, which he called {\em irreducible Levi graphs}.

In \cite{ALS1} Abreu, Labbate and Sheehan characterised {\em irreducible} pseudo $2$--factor isomorphic cubic bipartite graphs (and hence
gave a further partial answer to Conjecture \ref{quest1})
as follows:

    \begin{theorem}\label{MainThm}\cite{ALS1}
The Heawood and the Pappus graphs are the only irreducible
Levi graphs which are pseudo $2$--factor isomorphic.
    \end{theorem}

This approach is not feasible to prove Conjecture \ref{quest1} and hence the main Conjecture \ref{conj1} by studying
the $2$--factors of reducible configurations from the
set of $2$--factors of their underlying irreducible ones as the following discussion shows.

It is well known that the $7_3$ configuration, whose Levi graph is the Heawood graph, is not
Martinetti extendible and that the Pappus configuration is Martinetti extendible in a unique
way; it is easy to show that this extension is not pseudo-2-factor isomorphic.
Let $\cal C$ be a symmetric configuration $n_3$ and $\cal C$ be a symmetric configuration $(n + 1)_3$
obtained from $\cal C$ through a Martinetti extension. It can be easily checked that there are $2$--factors in $\cal C$ that cannot be reduced to a $2$--factor in $\cal C$. On the other hand, all of its Martinetti reductions are no longer pseudo $2$--factor isomorphic (for further details cf. \cite{ALS1}).

In the next section, we will see that Conjecture \ref{quest1} has been disproved while Conjecture \ref{conj1} still holds.

\section{A counterexample to the pseudo 2--factor's conjecture}\label{counterpseudo}

In this section we present the counterexample by Jan Goedgebeur to the pseudo $2$--factor isomorphic bipartite's Conjecture \ref{quest1} obtained using exhaustive research via parallel computers (for details refer to \cite{Go15}). Recently, in \cite{config30}, it was shown how it could be constructed from the Heawood graph and the generalised Petersen graph $GP(8, 3)$, which are the Levi graphs of the Fano $7_3$ configuration and the M\"obius-Kantor $8_3$ configuration, respectively. 

Using the program {\em minibaum} \cite{Br}, J.Goedgebeur generated all cubic bipartite graphs with girth at least $6$ up to $40$ vertices and all cubic bipartite graphs with girth at least $8$ up to $48$ vertices. The counts of these graphs can be found in \cite[Table 1]{Go15}. Some of these graphs can be downloaded from the House of Graphs Database \cite{hog2}, available online at \url{https://houseofgraphs.org/}.
He, then, implemented a program that tests if a given graph is pseudo 2--factor isomorphic and applied it to the generated cubic bipartite graphs. This yielded the following results:

\begin{remark}\cite{Go15}
There is exactly one essentially $4$--edge--connected pseudo $2$--factor isomorphic graph different from the Heawood
graph and the Pappus graph among the cubic bipartite graphs with girth at least $6$ with at most $40$ vertices.
\end{remark}

\begin{remark}\cite{Go15}
There is no essentially $4$--edge--connected pseudo $2$--factor isomorphic graph among the cubic bipartite graphs
with girth at least $8$ with at most $48$ vertices.
\end{remark}

This implies that Conjecture \ref{quest1} (and consequently also Conjecture \ref{quest2}) is false. However, since all $2$--factor Hamiltonian graphs are pseudo $2$--factor isomorphic and $\cal G$ is not $2$--factor Hamiltonian, this implies the following remark:

\begin{remark}\cite{Go15}
Conjecture \ref{conj1} holds up to at least $40$ vertices and holds for cubic bipartite graphs with girth at least $8$ up to at
least $48$ vertices.
\end{remark}

\

The counterexample found has $30$ vertices and there are no additional counterexamples up to at least $40$ vertices and also no counterexamples among the cubic bipartite graphs with girth at least $8$ up to at least $48$ vertices.

The counterexample $\cal G$ is stored in the House of Graphs Database \cite{hog2} by searching for the keywords \texttt{pseudo 2-factor isomorphic *counterexample} where it can be downloaded and several of its invariants can be inspected (see Figure \ref{fig:jan} below).

\
\begin{figure}
\centering
\begin{subfigure}{0.5\textwidth}
  \centering
  \includegraphics[width=.8\linewidth]{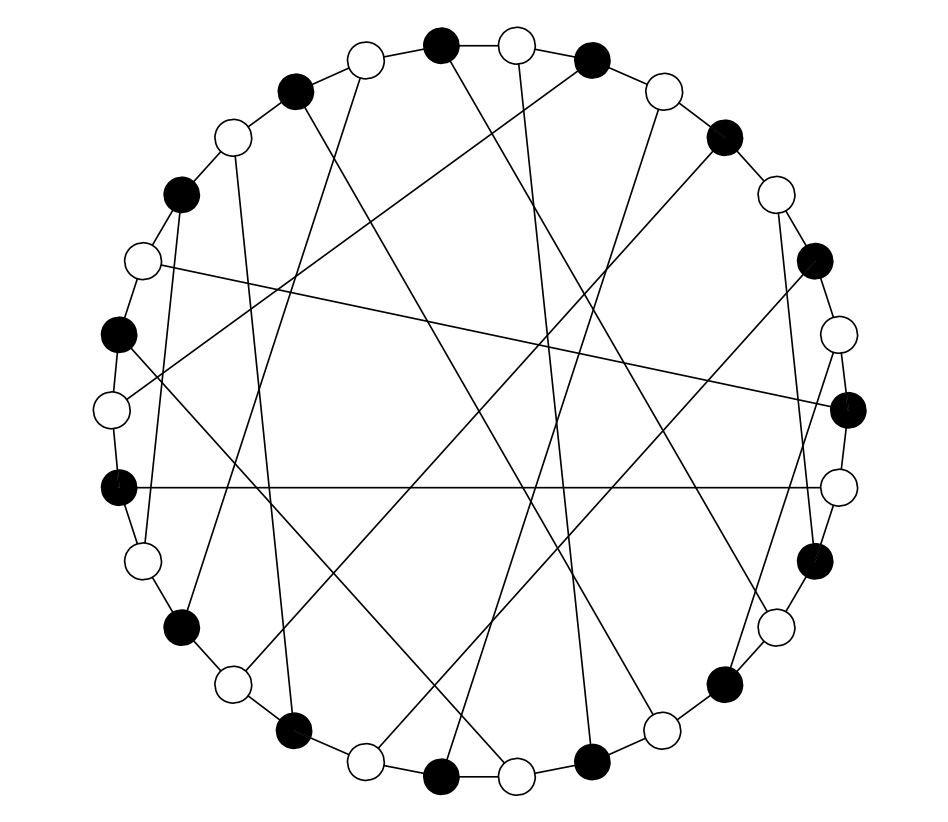}
  \caption{The graph $\cal{G}$, as appear in the House of Graphs Database \cite{hog2}.}
  \label{fig:jan}
\end{subfigure}%
\begin{subfigure}{0.5\textwidth}
  \centering
  \includegraphics[width=.9\linewidth]{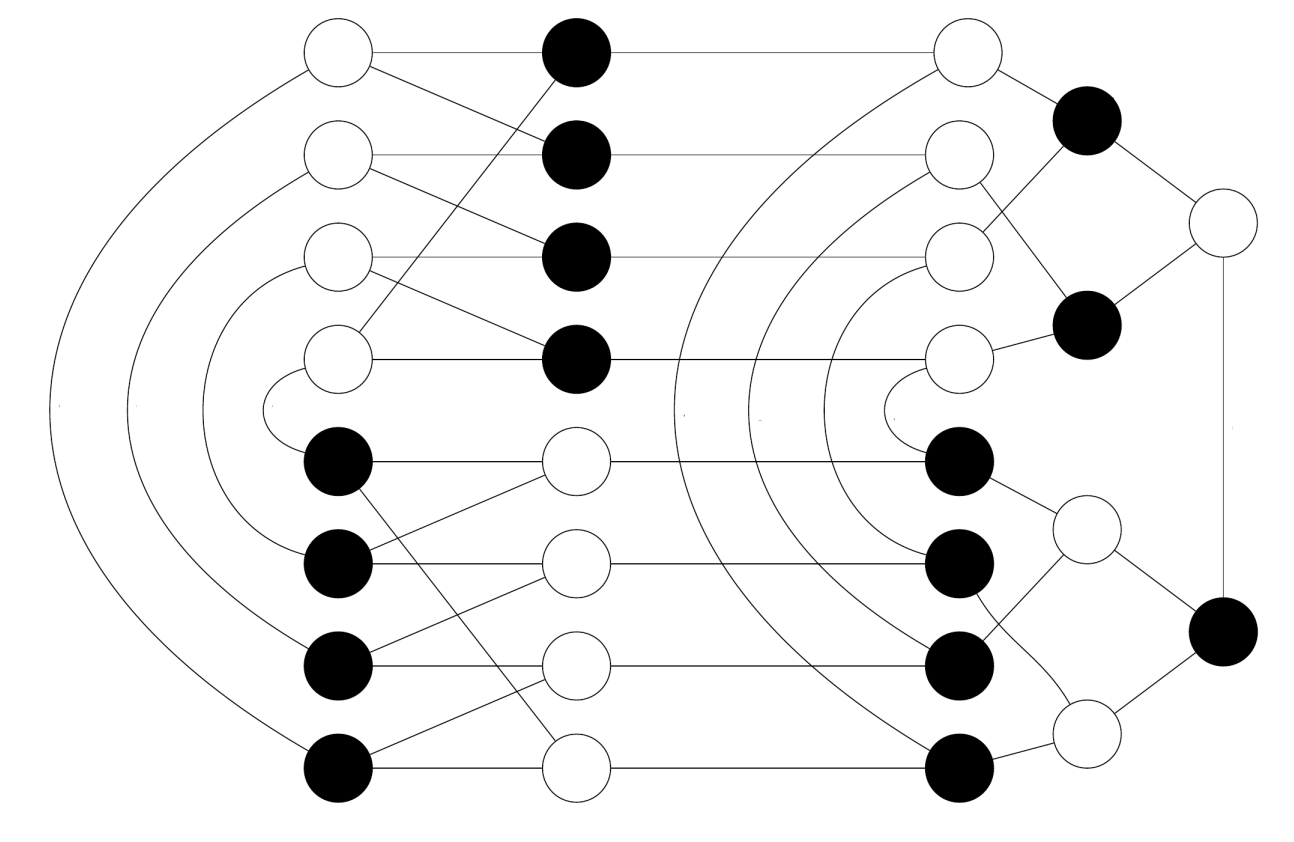}
  \caption{The geometric construction of $\cal{G}$, as drawn in \cite{config30}.}
  \label{fig:c30}
\end{subfigure}
\caption{Two possible drawings of the graph $\cal{G}$.}
\label{fig:counter}
\end{figure}

In what follows we will briefly describe the geometric construction of $\cal G$ presented in \cite{config30}. Consider the classical representation of the M\"obius-Kantor configuration as two quadrilaterals simultaneously inscribed and circumscribed (cf. \cite[pag. 430]{cox}). Disregarding the circumscription, i.e. removing the corresponding incidences it defines, we obtain an M\"obius-Kantor residue in
which the valency of 4 points and 4 lines decreases from three to two. Similarly, removing the incidences of a quadrilateral in the Fano configuration, we obtain a Fano residue with 4 points and 4 lines of valency two. The configuration $\cal C$  which has $\cal G$ as
its Levi Graph, then arises by adding incidences among points and lines of valency two between the Fano residue and the M\"obius-Kantor residue in a precise way (see Figure \ref{fig:c30}), fully described in \cite[Section 2]{config30}.

In detail, $\cal{G}$ has girth $6$, cyclic edge--connectivity $6$, it is not vertex--transitive, it has $312$ $2$--factors and the cycle sizes of its $2$--factors are: $(6,6,18)$, $(6,10,14)$, $(10,10,10)$ and $(30)$. Moreover, its automorphism group of order $144$ is isomorphic to $(\mathbb{Z}_3 \times \mathbb{Z}_3)\rtimes(D_4 \times \mathbb{Z}_2)$. The authors of \cite{config30} pointed out that the construction of joining residues of Levi graphs of $n_3$ configurations, does not preserve, in general a strong property such as being pseudo 2-factor isomorphic. Intuitively, there is too much room to produce 2-factors of both parities; whereas, the M\"obius-Kantor residue and the Fano residue are very compact. Moreover, also using several copies of the same residues does not preserve the behavior of the parity of cycles in a 2–factor.\\

\section{Strongly pseudo 2--factor isomorphic graphs}

Abreu, Labbate and Sheehan in \cite{ALS} have extended the above-mentioned results on regular pseudo 2--factor isomorphic bipartite graphs to the not necessarily bipartite case introducing the family of strongly pseudo 2--factor isomorphic graphs:

\begin{definition}\label{pseudo-strongly}
{\em Let $G$ be a graph which has a $2$--factor. For each $2$--factor
$F$ of $G$, let $t^*_i(F)$ be the number of cycles of $F$ of
length $2i$ modulo $4$. Set $t_i$ to be the function defined on the
set of $2$--factors $F$ of $G$ by:}

$$
t_i(F) = \left\{
\begin{array}{ll}
0 & \text{if $\, \, t^*_i(F)$ is even}\\
1 & \text{if $\, \, t^*_i(F)$ is odd}
\end{array}
\right. \quad \quad (i=0,1).
$$

\noindent {\em Then $G$ is said to be}
strongly pseudo $2$--factor isomorphic
{\em if both $t_0$ and $t_1$ are constant functions.
Moreover, if in addition $t_0=t_1$, set $t(G):=t_i(F)$, $i=0,1$.}
    \end{definition}

By definition, {\em if $G$ is strongly pseudo $2$--factor
isomorphic then $G$ is pseudo $2$--factor isomorphic.} On the other
hand there exist graphs such as the Dodecahedron which are pseudo
$2$--factor isomorphic but not strongly pseudo $2$--factor
isomorphic: the $2$--factors of the Dodecahedron consist either of
a cycle of length $20$ or of three cycles: one of length $10$
and the other two of length $5$.

In the bipartite case, pseudo $2$--factor isomorphic and
strongly pseudo $2$--factor isomorphic are equivalent.

In what follows we will denote by $HU$, $U$, $SPU$ and $PU$ the sets of
$2$--factor Hamiltonian, $2$--factor isomorphic, strongly pseudo $2$--factor isomorphic
and pseudo $2$--factor isomorphic graphs, respectively.
Similarly, $HU(k)$, $U(k)$, $SPU(k)$, $PU(k)$ respectively denote the
$k$--regular graphs in $HU$, $U$, $SPU$ and $PU.$

    \begin{theorem}\label{di-degree}\cite{ALS}
Let $D$ be a digraph with $n$ vertices and $X$ be a directed
$2$--factor of $D$. Suppose that either
\begin{enumerate}[(a)]
  \item $d^+(v) \geq \lfloor log_2 \; n \rfloor+2$ for all $v \in V(D)$, or
  \item $d^+(v)=d^-(v) \geq 4$ for all $v \in V(D)$
\end{enumerate}

\noindent Then $D$ has a directed $2$--factor $Y$ with a different
parity of number of cycles from $X$.
    \end{theorem}

Let $DSPU$ and $DPU$ be the sets of digraphs in $SPU$
and $PU$, i.e. strongly pseudo and pseudo $2$--factor isomorphic digraphs, respectively.
Similarly, $DSPU(k)$ and $DPU(k)$ respectively denote the $k$--diregular digraphs
in $DSPU$ and $DPU$.

    \begin{cor}\label{DSPU-DQUempty}\cite{ALS}
\begin{enumerate}[(i)]
  \item $DSPU(k)=DPU(k)= \emptyset$ for $k \geq 4;$
  \item If $D \in DPU$ then $D$ has a vertex
  of out--degree at most $\lfloor log_2 \; n \rfloor + 1.$
\end{enumerate}
\end{cor}

    \begin{theorem}\label{G-degree}\cite{ALS}
Let $G$ be a graph with $n$ vertices and $X$ be a
$2$--factor of $G$. Suppose that either
\begin{enumerate}[(a)]
  \item $d(v) \geq 2(\lfloor log_2 \; n \rfloor+2)$ for all $v \in V(G)$, or
  \item $G$ is a $2k$--regular graph for some $k \geq 4$
\end{enumerate}

\noindent Then $G$ has a $2$--factor $Y$ with a different
parity of number of cycles from $X$.
    \end{theorem}

 \begin{cor}\label{QUempty}\cite{ALS}
\begin{enumerate}[(i)]
  \item If $G \in PU$ then $G$ contains a vertex of degree at most
  $2\lfloor log_2 \; n \rfloor + 3;$
  \item $PU(2k)=SPU(2k)=\emptyset$ for $k \geq 4.$
\end{enumerate}
\end{cor}

We know that there are examples of graphs in $PU(3)$, $SPU(3)$, $PU(4)$ and $SPU(4)$, hence they are not empty
and we have seen (cf. Conjecture \ref{hukemptyset}) that it has been conjectured in \cite{AAFJLS} that $HU(4)=\{K_5\}$.

\noindent There are many gaps in our knowledge even when we restrict attention to regular graphs. Some questions arise naturally. A few of them are listed below.

    \begin{prob}\label{emptyodd}
Is $PU(2k+1)=\emptyset$ for $k\ge 2$?
    \end{prob}

In particular, it is wondering if  $PU(7)$ and $PU(5)$ are empty.

\begin{prob}\label{furthprob1} Is $PU(6)$ empty? \end{prob}

\begin{prob}\label{furthprob2} Is $K_5$ the only $4$--edge--connected graph in $PU(4)$? \end{prob}

In \cite{ALS}, relations between pseudo strongly 2--factor isomorphic graphs and a class of graphs called {\em odd 2--factored snarks} are investigated. The next section is devoted to this class of snarks.

\section{Odd 2--factored snarks}\label{odd2fact}

A {\em snark} (cf. e.g.~\cite{HS}) is a bridgeless cubic graph with chromatic
index four (by Vizing's theorem the chromatic index of every cubic
graph is either three or four so a snark corresponds to the special case
of four). In order to avoid trivial cases, snarks are usually assumed to have
girth at least five and not to contain a non--trivial $3$--edge cut (i.e. they are cyclically $4$--edge connected).

Snarks were named after the mysterious and elusive creature in Lewis Caroll's famous poem {\em The Hunting of The Snark} by Martin Gardner in $1976$
\cite{G76}, but it was P. G. Tait in 1880 that initiated the study of snarks, when he proved that the four colour theorem is equivalent to the statement that {\em no snark is planar}~\cite{Ta1880}. The Petersen graph $P$ is the smallest snark and Tutte conjectured that all snarks have Petersen graph minors.
This conjecture was proven by Robertson, Seymour and Thomas (cf. \cite{RST2}).
Necessarily, snarks are non--Hamiltonian.

The importance of the snarks does not only depend on the four colour theorem. Indeed, there are several important open problems
such as the classical cycle double cover conjecture~\cite{Sey79,Sze73}, Fulkerson's conjecture~\cite{Fu71} and Tutte's
5--flow conjecture~\cite{Tu54} for which it is sufficient to prove them for snarks.
Thus, minimal counterexamples to these and other problems must reside, if they exist at all, among the family of snarks.

At present, there is no uniform theoretical method for studying snarks and their behaviour. In particular, little is known about the structure of $2$--factors in a given snark.

Snarks play also an important role in characterising regular graphs with some conditions imposed on their 2--factors.
Recall that a $2$--factor is a $2$--regular spanning subgraph of a graph $G$.

We say that a graph $G$ is {\em odd $2$--factored} (cf.~\cite{ALS}) if for each $2$--factor $F$ of $G$ each cycle of $F$ is odd.

By definition, {\em an odd $2$--factored graph $G$ is pseudo $2$--factor isomorphic}.
Note that,  odd $2$--factoredness is not the same as the {\em oddness} of a (cubic) graph (cf. e.g.\cite{Z}).

\begin{lemma}\cite{ALS}\label{snark1}
Let $G$ be a cubic $3$--connected odd $2$--factored  graph then $G$ is a snark.
\end{lemma}

In~\cite{ALS}, it was studied which snarks are odd $2$--factored and the following conjecture has been posed:

\begin{conj}\cite{ALS}\label{con1}
A snark is odd $2$--factored if and only if $G$ is the
Petersen graph, Blanu\v{s}a~$2$, or a Flower snark $J(t)$, with $t \ge 5$ and odd.
\end{conj}

In \cite{ALRS}, it was shown a general construction, which we report below, of odd $2$--factored snarks performing the Isaacs' dot--product~\cite{I75} on edges with particular properties, called {\em bold--edges} and {\em gadget--pairs} respectively, of two snarks $L$ and $R$.

\

\noindent{\sc Construction: Bold--Gadget Dot Product.}\cite{ALRS}

\

\begin{itemize}
\item Take two snarks $L$ and $R$ with bold--edges (cf. definitionition~\ref{boldedge}) and gadget--pairs (cf. definitionition~\ref{gadgetpair}), respectively;
\item Choose a bold--edge $xy$ in $L$;
\item Choose a gadget--pair $f$, $g$ in $R$;
\item Perform a dot product $L \cdot R$ using these edges;
\item Obtain a new odd 2--factored snark (cf. Theorem~\ref{OddDot}).
\end{itemize}

Note that in what follows the existence of a 2--factor in a snark is guaranteed since they are bridgeless by definition.

\begin{definition}\label{boldedge}\cite{ALRS}
Let $L$ be a snark. A {\em bold--edge} is an edge $e=xy \in L$ such that the following conditions hold:

\begin{enumerate}[(i)]
  \item All $2$--factors of $L-x$ and of $L-y$ are odd;
  \item all $2$--factors of $L$ containing $xy$ are odd;
  \item all $2$--factors of $L$ avoiding $xy$ are odd.
\end{enumerate}
\end{definition}

Note that not all snarks contain bold--edges (cf. \cite[Proposition 4.2]{ALRS} and \cite[Lemma 5.1]{ALRS}).
Furthermore, conditions $(ii)$ and $(iii)$ are trivially satisfied if $L$ is odd $2$--factored.

\begin{lemma}\label{PetBold}\cite{ALRS}
The edges of the Petersen graph $P_{10}$ are all bold--edges.
\end{lemma}

\begin{definition}\label{gadgetpair}\cite{ALRS}
Let $R$ be a snark. A pair of independent edges $f=ab$ and $g=cd$ is called a {\em gadget--pair} if the following conditions hold:

\begin{enumerate}[(i)]
  \item There are no $2$--factors of $R$ avoiding both $f,g$;
  \item all $2$--factors of $R$ containing exactly one element of $\{f,g\}$ are odd;

  \item all $2$--factors of $R$ containing both $f$ and $g$ are odd. Moreover, $f$ and $g$ belong to different cycles in each such factor.

  \item all $2$--factors of $(R-\{f,g\}) \cup \{ac,ad,bc,bd\}$ containing exactly one
  element of $\{ac,ad,bc,bd\}$, are such that the cycle containing the new edge is even and all other cycles are odd.
\end{enumerate}
\end{definition}

Note that, finding gadget--pairs in a snark
is not an easy task and, in general, not all snarks contain gadget--pairs (cf. \cite[Lemma 5.2]{ALRS}).

\begin{figure}[h]
\begin{center}
\includegraphics[scale=0.3]{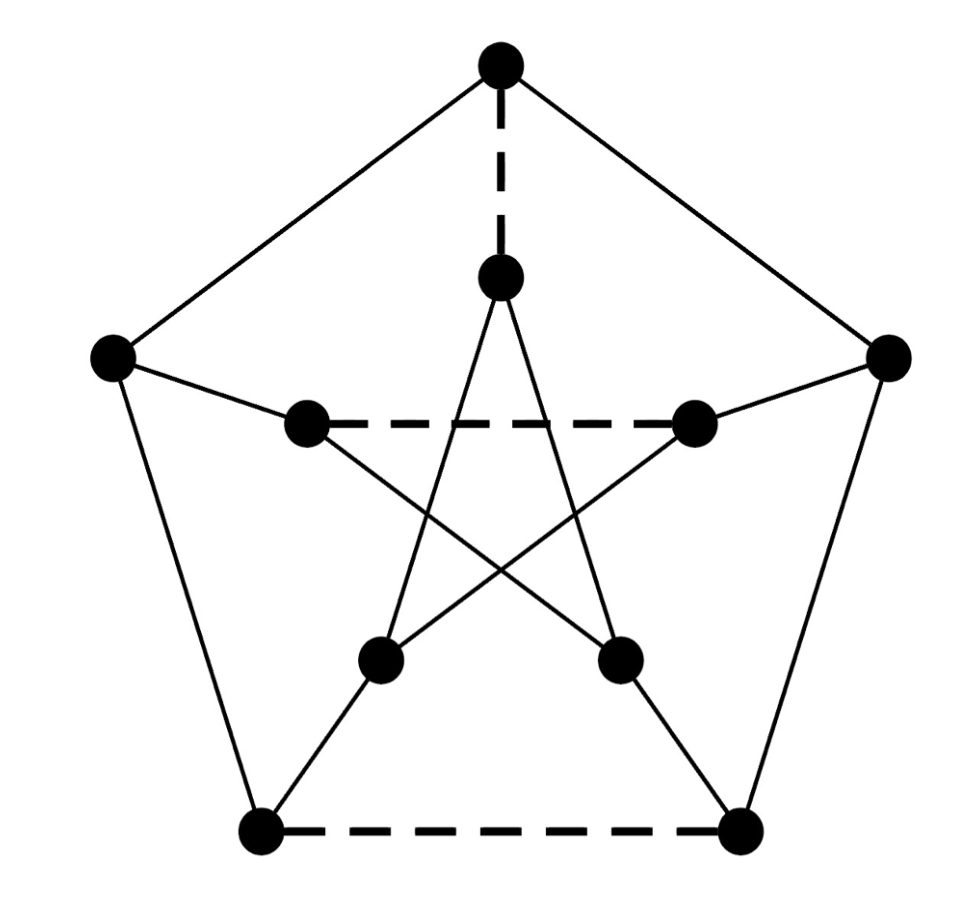}
\caption{Any pair of the dashed edges is a gadget--pair in $P_{10}$}
\label{PetSpc}
\end{center}
\end{figure}

Let $H:=\{x_1y_1,x_2y_2,x_3y_3\}$ be the two horizontal
edges and the vertical edge respectively (in the pentagon--pentagram representation) of $P_{10}$ (cf. Figure~\ref{PetSpc}).

\begin{lemma}\label{PetGadg}\cite{ALRS}
Any pair of distinct edges $f,g$ in the set $H$ of $P_{10}$ is a gadget--pair.
\end{lemma}

The following theorem is used to construct new odd $2$--factored snarks.

\begin{theorem} \label{OddDot}\cite{ALRS}
Let $xy$ be a bold--edge in a snark $L$ and let $\{ab,cd\}$ be a gadget--pair in a snark $R$.
Then $L \cdot R$ is an odd $2$--factored snark.
\end{theorem}

In particular, without going into lengthy details (the interested reader might find those in \cite{ALRS}), this method allowed to construct two instances of odd $2$--factored snarks of order
$26$ and $34$ isomorphic to those obtained by Brinkmann et al. in \cite{BGHM} through an exhaustive computer search on all snarks of order $\le 36$ that has allowed them to disprove the above conjecture (cf. Conjecture~\ref{con1}). 


To approach the problem of characterising all odd $2$--factored snarks, it was considered
the possibility of constructing further odd $2$--factored snarks with the technique presented above, which
relies on finding other snarks with bold--edges and/or gadget--pairs,
The results obtained so far give rise to the following partial characterization:

\begin{theorem}\label{partialodd}\cite{ALRS}
Let $G$ be an odd $2$--factored snark of cyclic edge--connectivity four that can be constructed from the Petersen graph and the Flower snarks using the bold--gadget dot product construction. Then $G \in \{P_{18}$, $P_{26}$, $P_{34}\}$.
\end{theorem}

Finally, in \cite{ALRS} a new conjecture about odd $2$--factored snarks has been posed.

\begin{conj}\label{newconj}\cite{ALRS}
Let $G$ be a cyclically $5$--edge connected odd $2$--factored snark. Then $G$ is either the Petersen graph or the Flower snark $J(t)$, for odd $t\ge5$.
\end{conj}

\begin{remark}\label{remconj}

\noindent $(i)$ A minimal counterexample to Conjecture~\ref{newconj} must be a cyclically $5$--edge connected snark
of order at least $36$.
Moreover, as highlighted in~\cite{BGHM}, order $34$ is a turning point for several properties of snarks.

\noindent $(ii)$ It is very likely that, if such counterexample exists, it will arise from the superposition operation by M.Kochol \cite{K} applied to one of the known odd $2$--factored snarks.

\noindent $(iii)$ J. Goedgebeur \cite{Go} checked that none of the snarks (in particular those with girth $6$ of order $38$) that G.Brinkmann and himself generates in \cite{BG} is an odd $2$--factored snark.


\end{remark}


\end{document}